\documentclass[french]{amsart}

\usepackage[french]{babel}
\usepackage[applemac]{inputenc}
\usepackage{latexsym}

\usepackage{amsmath}
\usepackage{amsfonts}
\usepackage{amssymb}
\usepackage{verbatim}
\usepackage{psfrag}
\usepackage{graphicx}

\theoremstyle{theorem}
\newtheorem{lemm}{Lemme}[section]
\newtheorem{prop}[lemm]{Proposition}

\newtheorem{theo}[lemm]{Th\'eor\`eme}
\theoremstyle{remark}
\newtheorem{defi}[lemm]{D\'efinition}

\newtheorem{rema}[lemm]{Remarque}

  \usepackage[all]{xy}
\xyoption{matrix}

\newcommand{\PGL}{\mathrm{PGL}}
\newcommand{\PSL}{\mathrm{PSL}}

\newcommand{\Bir}{\mathrm{Bir}}
\newcommand{\Cr}{\mathrm{Cr}}

\newcommand{\Aut}{\mathrm{Aut}}
\newcommand{\p}{\mathbb{P}}

\newcommand{\GL}{\mathrm{GL}}

\newcommand{\J}{\mathit{J}}

\title{Groupes de Cremona, connexit\'e et simplicit\'e}
\date{\today}
\author{J\'er\'emy Blanc}

\begin{document}
\begin{abstract}Le groupe de Cremona est connexe en toute dimension et, muni de sa topologie, il est simple en dimension $2$.

\medskip

\noindent {\tiny\sc ABSTRACT} The Cremona group is connected in any dimension and, endowed with its topology, it is simple in dimension $2$.
\end{abstract}

\maketitle\begin{center}{\today}\end{center}
\section{Questions et résultats}
Soit $k$ un corps algébriquement clos. On note $\Cr_{n}(k)$ le groupe de Cremona de dimension $n$, groupe des transformations birationnelles de $\mathbb{P}^n_{k}$, anti-isomorphe à $\Aut_{k}(k(x_{1},\dots,x_{n}))$.
Ce groupe est muni d'une topologie naturelle (décrite à la section~\ref{Sec:Intro2}).

En 1974, dans un rapport  sur les questions ouvertes importantes en géométrie algébrique \cite{bib:Mum}, D.~Mumford consacre un paragraphe au groupe $\Cr_{2}(k)$. Il parle de mettre une topologie sur le groupe, et pose alors la question: ce groupe est-il simple? Le théorème~\ref{Thm:Cr2Simple} démontrée plus bas permet de répondre par l'affirmative.

La technique utilisée pour cela est élémentaire. Elle permet également de prouver que le groupe $\Cr_{n}(k)$ est connexe pour tout $n$ (Théorème~\ref{Thm:Connexite}). Ceci répond à une question posée par J.-P. Serre lors du 1000ème exposé Bourbaki
\cite{bib:SerreBour}, concernant la dimension $n\ge 3$, le cas 
$n\le 2$ étant déjà bien connu.

\smallskip

Cet article est articulé ainsi: la section~\ref{Sec:Intro2} donne des rappels sur la topologie de Zariski de $\Cr_{n}(k)$, la section~\ref{Sec:Prelim} présente un lemme de déformation, qui permet de montrer la simplicité de $\Cr_{2}(k)$ (section~\ref{Sec:Simpl}) et la connexité de $\Cr_{n}(k)$ (section~\ref{Sec:Connex}).

\medskip

Je tiens à remercier J.-P. Furter pour  des discussions intéressantes sur cet article, et tout spécialement J.-P. Serre pour ses relectures attentives de cet article et ses précieuses corrections.
\section{La topologie de Zariski de $\Cr_{n}(k)$}\label{Sec:Intro2}
Soit $X$ une $k$-variété ($k$ est toujours le corps algébriquement clos fixé au départ). On note $\Bir(X)$ l'ensemble des applications birationnelles $X\dasharrow X$, et $\Aut(X)\subset \Bir(X)$ le groupe des automorphismes de $X$. 

\smallskip

Afin de d\'ecrire la topologie de $\Bir(X)$, d\'ecrivons tout d'abord les  morphismes  $A\to \Bir(X)$: 
\begin{defi}\label{defi:Famille}
Une \emph{famille algébrique} d'éléments de $\Bir(X)$ est la donnée d'une application rationnelle  $f:A\times X\dasharrow X$ où $A$ est une $k$-variété, définie sur un ouvert dense $U$ tel que pour tout $a\in A$, $U_{a}:=U\cap (\{a\}\times A)$ soit un ouvert dense de $\{a\}\times A$ et la restriction de $f$ \`a cet ouvert soit un isomorphisme de $U_a$ sur un ouvert dense de $X$.

Pour tout $a\in A$, l'application birationnelle $x\dasharrow f(a,x)$ représente alors un élément $f_{a}\in \Bir(X)$.
La famille $f_a$ $(a\in A)$ représente une application $A\to \Bir(X)$, que l'on appelera \emph{morphisme} de $A$ vers $\Bir(X)$.
\end{defi}
Cette d\'efinition correspond \`a celle de \cite{bib:SerreBour} et \cite[\S 1]{bib:Demazure}; un morphisme  $A\to \Bir(X)$ correspond alors à un pseudo-automorphisme du $A$-schéma $A\times X$.
On d\'efinit la topologie de Zariski sur $\Bir(X)$ de la manière suivante (voir \cite[\S 1.6]{bib:SerreBour}):

\begin{defi}
On dit qu'un ensemble $R\subset\Bir(X)$ est \emph{fermé} si pour toute $k$-variété $A$ et tout morphisme $A\to\Bir(X)$, la préimage de $R$ est fermée.
\end{defi}
Comme l'explique \cite{bib:SerreBour}, ceci donne une topologie sur $\Bir(X)$, qui est la topologie la plus fine qui rende les morphismes vers $\Bir(X)$ continus. De plus, en définissant de manière analogue la topologie de Zariski sur $\Bir(X)\times \Bir(X)$, la composition donne une application continue $\Bir(X)\times \Bir(X)\to \Bir(X)$. 

En particulier, on peut restreindre ceci à $\Aut(X)$ et mettre ainsi une topologie sur ce groupe. Lorsque $X=\mathbb{P}^n_{k}$, on peut démontrer que l'on retrouve la topologie de Zariski habituelle du groupe algébrique $\Aut(X)=\PGL(n+1,k)$ et qu'en fait $\Aut(X)\to \Cr_n(X)$ est une immersion fermée \cite{bib:SerreLettre}.



Les groupes qui nous intéressent le plus sont ceux où la $k$-variété $X$ est rationnelle. On rappelle que si $X$ est rationnelle, de dimension $n$, alors $\Bir(X)$ s'identifie naturellement à $\Cr_n(k)=\Bir(\mathbb{P}^n_{k})$ via une application birationnelle choisie $X\dasharrow \mathbb{P}^n_k$; le choix de celle-ci fait juste varier l'homéomorphisme $\Bir(X)\to \Cr_n(k)$.
Dans la suite, on prendra le plus souvent $X=\mathbb{A}^n_{k}$ ou $X=\mathbb{P}^n_{k}$, suivant les besoins.

\section{Préliminaires techniques}\label{Sec:Prelim}
\subsection{Le groupe de de Jonquières}
Pour $n\geq 2$, notons $\phi$ la projection \[(x_0:\dots:x_n)\dasharrow (x_1:\dots:x_n)\] de $\mathbb{P}^n_k$ dans $\mathbb{P}^{n-1}_{k}$.
On  appelle \emph{groupe de de Jonqui\`eres} $\J_n$ le sous-groupe de $\Bir(\mathbb{P}^n_k)$ qui préserve l'ensemble des fibres de $\phi$. On note $\J_n^0$ le sous-groupe de $\J_n$ constitué des éléments qui préservent une fibre générale de $\phi$.

De manière affine, on peut restreindre $\phi$ à la projection  $k^{n}\to k^{n-1}$, et ainsi voir que $\J_n$ est naturellement isomorphe à $\J_n^0\rtimes\Bir(k^{n-1})$, où \[\J_n^0\simeq \Aut(\mathbb{P}^1_K)\simeq\PGL(2,K),\mbox{ avec }K=k(x_1,\dots,x_{n-1}).\]
L'homomorphisme déterminant $\GL(2,K)\to K^{*}$ induit un homomorphisme surjectif 
\[\det\colon \PGL(2,K)\to K^{*}/(K^{*})^2,\]
où $(K^{*})^2$ désigne l'ensemble des carrés de $K^{*}$.
On notera $\J^1_n\subset \J^0_n$ le sous-groupe normal correspondant au noyau de $\det$. Alors, l'homomorphisme précédent nous donne \[\J^1_n\simeq \PSL(2,K).\]

 Le groupe $\J^1_n$ est simple \cite[Chapitre II, $\S 2$]{bib:Dieudonne}. De plus, comme tout élément $f\in \J^0_n$ satisfait $\det(f^2)=1$, le quotient $\J^0_n/\J^1_n$ est un groupe abélien de type $(2,\dots,2,\dots)$. Les classes de $\J^0_n \pmod{\J^1_n}$ sont représentées par les involutions de de Jonquières $f_h\colon (x_1,\dots,x_n) \dasharrow (x_1,\dots,x_{n-1},h/x_n)$, où $h\in k(x_1,\dots,x_{n-1})^{*}=K^{*}$. Puisque $\det(f_h)=-h$; deux involutions $f_h$ et $f_{h'}$ représentent la même classe si et seulement si $h/h'$ est un carré dans $K^{*}$.

\subsection{Dérivée normale}
\label{SubSec:Deriveenormale}
Dans cette section, on se donne la situation suivante:

Partons d'une $k$-variété lisse $X$. Soit  $Y$ la droite affine sur $k$ et soit $Z=X\times Y$; le morphisme $x\mapsto (x,0)$ identifie $X$ à une sous-variété de $Z$. Soit $U$ un ouvert de $Z$ et soit $f\colon U\to Z$ un morphisme qui applique $X$ dans lui-même; notons $f_X\colon U\to X$ et $f_Y\colon U\to Y$ le deux composantes de $f$. 

À partir de cette donnée, on va définir la dérivée normale de $f$, qui est un morphisme $f_0\colon Z\to Z$, et montrer qu'il s'agit d'une limite de conjugués de $f$.

\medskip

La fonction $f_Y$ a la propriété que $f_Y(x,y)=0$ si $y=0$; on en déduit que $f_Y$ est divisible par la fonction "y", ce qui veut dire que $f_Y(x,y)=y\cdot g_Y(x,y)$ pour une certaine fonction $g_Y$ sur $U$. Ceci nous permet de définir la \emph{dérivée normale} de $f$, qui est le morphisme
\[f_0\colon Z\to Z,\]
donné par la formule $f_0(x,y)=(f_X(x,0),y\cdot g_Y(x,0))$.

\medskip

On remarque que $f_0$ ne dépend que du comportement de $f$ dans un voisinage infinitésimal de $X$ et est une sorte de linéarisation de $f$. De plus $(x,0)\mapsto f_X(x,0)$ est la restriction de $f$ à $X$, ce qui implique que $f_0$ est compatible avec la projection $Z\to X$; on a le diagramme commutatif
\begin{equation}\label{DiagComm1}\xymatrix{Z\ar[d] \ar@{-->}[r]^{f_0}& Z\ar[d]\\ X\ar@{-->}[r]^{f|_X}&X. }\end{equation}
De plus, $f_0$ est une homothétie sur chaque fibre.

\medskip

Montrons maintenant que $f_0$ est une limite des conjugués de $f$.  Soit $T$ un autre exemplaire de la droite affine sur $k$; pour tout $t\in T$, soit $U_t$ l'ouvert de $Z$ formé des $(x,y)$ tels que $(x,ty)$ appartienne à $U$. La réunion $U_T$ des $(t,U_t)$ est un ouvert de $T\times Z$ contenant $T\times X$. Si $t\not=0$, soit $s_t$ l'automorphisme $(x,y)\mapsto (x,ty)$ de $Z$ et posons $f_t=s_t^{-1}\circ f \circ s_t$,  qui a un sens sur $U_t$; si $t=0$ définissons $f_t=f_0$ comme ci-dessus.
\begin{lemm}\label{lem:FUT}
Avec les notation précédentes, la famille des $f_t$ $(t\in T)$ définit un morphisme $F\colon U_T\to Z$.
\end{lemm}
\begin{proof}
On a $F(t,x,y)=(f_X(x,ty),y\cdot g_Y(x,ty))$: lorsque $t\not=0$ cela résulte de $t^{-1}f_Y(x,ty)=y\cdot g_Y(x,ty)$ et lorsque $t=0$ c'est la définition de $f_0$. Le lemme suit alors du fait que $(t,x,y)\mapsto f_X(x,ty)$ et $(t,x,y)\mapsto g_Y(x,ty)$ soient des morphismes définis sur $U_T$. 
\end{proof}
\begin{lemm}\label{Lem:DefoMorphismeBir}
Avec les mêmes notations qu'avant, supposons de plus que $X$ est irréductible et que $f$ soit un isomorphisme de $U$ sur un ouvert $V$ de $Z$ $($ce qui implique que $f$ est birationnelle$).$

Alors, la famille $f_t$ $(t\in T)$ définit un morphisme $T\to \Bir(Z)$ $($au sens de la définition~$\ref{defi:Famille})$.
\end{lemm}
\begin{proof}
Le lemme~\ref{lem:FUT} montre que $F\colon U_T\to Z$ est un morphisme, qui induit donc une application rationnelle $T\times Z\dasharrow Z$. Pour tout $t\in T$, $U_T\cap (\{t\}\times Z)$ n'est rien d'autre que $\{t\}\times U_t$, ouvert dense de $\{t\}\times Z$ par construction, et la restiction de $F$ à cet ouvert correspond à $f_t$. Il reste à voir que $f_t\colon U_t\to Z$ est un isomorphisme sur un ouvert dense de $Z$, pour tout $t\in T$. 

Notons $r\colon V\to U$ l'inverse de $f$, qui applique $X$ dans lui-même par construction, et utilisons la construction précédente pour $r=(r_X,r_Y)$. On a $V_t=\{(x,y)\in Z\ |\ (x,ty)\in V\}$ et le lemme~\ref{lem:FUT} nous donne un morphisme $R\colon V_T\to Z$, dont la restriction à $\{t\}\times V_t$ correspond à $r_t$.

Par construction, on sait que pour $t$ non-nul,  $r_t=s_t^{-1}\circ r \circ s_t$ est l'inverse de $f_t=s_t^{-1}\circ f \circ s_t$.

Le morphisme $V_T\to Z$ donné par $(t,x,y)\mapsto F(t,R(t,x,y))$  se restreint à l'identité de $\{t\}\times V_t$ pour tout $t$ non-nul. À la limite, ce morphisme est donc également l'identité pour $t=0$. En faisant de même pour $(t,x,y)\mapsto R(t,F(t,x,y))$, on voit que $f_t \circ r_t$ et $r_t\circ f_t$ sont l'identité sur respectivement $V_t$ et $U_t$ pour tout $t\in T$, ce qui achève la démonstration.

 Le lecteur peut également remarquer que ces deux relations peuvent se déduire de la forme explicite de $F(t,x,y)$ et $R(t,x,y)$ donnée dans la preuve du lemme~\ref{lem:FUT}.
\end{proof}

\subsection{Le lemme de déformation appliqué au groupe de Cremona}

Rappelons que si $Z$ est une variété irréductible lisse, si $f\in \Bir(Z)$ et $H,H'\subset Z$ sont deux hypersurfaces irréductibles, on dit que $f$ \emph{se restreint à une application birationnelle $f|_{H}\colon H\dasharrow H'$} si $f$ est définie sur un ouvert $U$ tel que $U\cap H$ soit un ouvert dense de $H$ et tel que $f|_{U\cap H}\colon U\cap H\to H'$ soit une immersion ouverte. On peut également présenter cette notion de la façon suivante: les hypersurfaces $H$ et $H'$ définissent des valuations discrètes $v$ et $v'$ du corps des fonctions de $Z$, et l'on demande que $f$ transforme $v$ en $v'$.

En appliquant les résultats de la section $\ref{SubSec:Deriveenormale}$ au cas d'une transformation birationnelle de $\mathbb{P}^n_k$, on trouve le résultat suivant.

\begin{lemm}\label{Lem:DefoCremona}
Pour $n\geq 2$, notons $H_0\subset \mathbb{P}^{n}_k$ l'hyperplan d'équation $x_0=0$. Soit $f\in \Bir(\mathbb{P}^n_k)$ un élément qui se restreint à une application birationnelle $f|_{H_0}\colon H_0\dasharrow H_0$.

Notons $Z\subset \mathbb{P}^n_k$ le complémentaire de l'hyperplan d'équation $x_n=0$ et $X=Z\cap H_0$, de telle sorte que $Z=X\times Y$ avec $Y\cong \mathbb{A}^1_k$. On se donne $U,V\subset Z\subset \mathbb{P}^n$ deux ouverts tels que $f$ se restreigne à un isomorphisme $U\to V$.

Alors, en reprenant les notations de la section~$\ref{SubSec:Deriveenormale}$, la dérivée normale $f_0$ de $f$ est un élément de $\J_n$ tel que le diagramme suivant soit commutatif

\begin{equation}\label{DiagComm2}\xymatrix{\p^n_k\ar@{-->}[d] \ar@{-->}[r]^{f_0}& \p^n_k\ar@{-->}[d]\\ H_0\ar@{-->}[r]^{f|_{H_0}}&H_0, }\end{equation}
où les flèches verticales correspondent à la projection $(x_0:x_1:\dots:x_n)\dasharrow (0:x_1:\dots:x_n)$.
Le morphisme donné par le lemme~$\ref{Lem:DefoMorphismeBir}$
\[\begin{array}{rcl}\mathbb{A}^1_k\cong T&\to &\Bir(Z)\cong \Bir(\mathbb{P}^n_k)\\
t&\mapsto &f_t
\end{array}\] envoie $0$ sur $f_0$, $1$ sur $f_1=f$ et $t\not=0$ sur 
$f_t=s_t^{-1}\circ f \circ s_t$
où $s_t$ correspond ici à l'automorphisme 
$(x_0:x_1:\dots:x_n)\mapsto (tx_0:x_1:\dots: x_n)$
de $\mathbb{P}^n_k$.
\end{lemm}
\begin{proof}
On a $X\cong \mathbb{A}^{n-1}_k$ et $Z=X\times Y\cong \mathbb{A}^n_k$ est un ouvert dense de $\mathbb{P}^n_k$. Comme $f$ se restreint à un isomorphisme $U\to V$ qui applique $X$ dans lui-même, on peut utiliser tous les résultats de la section~\ref{SubSec:Deriveenormale}. La projection $Z\to X$ correspond à $\phi\colon \p^n_k\dasharrow H_0$ donné par $(x_0:x_1:\dots:x_n)\dasharrow (0:x_1:\dots:x_n)$. La commutativité du diagramme $(\ref{DiagComm1})$ entraîne donc celle de $(\ref{DiagComm2})$ et implique que $f_0\in \J_n$.  Le morphisme $t\mapsto f_t$ est donné par le lemme~\ref{Lem:DefoMorphismeBir} et sa description ici suit directement de celle donnée précédemment.
\end{proof}

\section{Simplicité de $\Bir(\mathbb{P}^2_{k})$}\label{Sec:Simpl}
\begin{prop}
Supposons $n\geq 2$. Soit $N\subset \Cr_{n}(k)$ un sous-groupe non-trivial qui soit à la fois normal et fermé.
Alors, $N$ contient $\Aut(\mathbb{P}^n_k)\simeq \PGL(n+1,k)$ et $\J^1_{n}\simeq \PSL(2,k(x_1,\dots,x_{n-1}))$.
\end{prop}
\begin{proof}
On se donne un élément non-trivial $h\in N$, qui se restreint à un isomorphisme $h|_U:U\to V$, où $U,V$ sont des ouverts denses de $\mathbb{P}^n_k$. Soit $p$ un point de $U$ et notons $q=h(p)\in V$ son image; on peut supposer que $q$ et $p$ sont différents. Il existe un élément $\alpha \in \Aut(\mathbb{P}^n_k)$ qui fixe à la fois $q$ et $p$. Par conséquent, $g=(\alpha h^{-1}\alpha^{-1})h\in N$ fixe $p$ (et $q$). 

Notons $T_p$ l'espace tangent à $p$ et $\mathbb{P}(T_p)\simeq \p^{n-1}_k$ son projectivisé. Alors, $g$ induit un automorphisme $g_p\in \Aut(\mathbb{P}(T_p))$.
Montrons maintenant que pour un choix convenable de $\alpha$, l'automorphisme $g_p$ est non trivial. Comme $h$ envoie $p$ sur $q$ via un isomorphisme local, il induit un isomorphisme (linéaire) $l:\mathbb{P}(T_p)\to \mathbb{P}(T_q)$. En notant $\alpha_p\in \Aut(\mathbb{P}(T_p))$ et $\alpha_q\in \Aut(\mathbb{P}(T_q))$ les automorphismes induits par les actions respectives de $\alpha$ sur $\mathbb{P}(T_p)$ et $\mathbb{P}(T_q)$, on a $g_p=\alpha_p l^{-1}(\alpha_q)^{-1} l$. Pour que $g_p$ soit non trivial, il suffit par exemple de choisir $\alpha=(x_0:x_1:\dots:x_n)\mapsto (\lambda x_0: x_1:\dots :x_n)$, avec $\lambda\in k\backslash \{0,1\}$, si $p=(1:0:\dots:0)$ et $q=(0:1:0:\dots:0)$.

Soit $\sigma\colon (x_0:\dots:x_n)\dasharrow (\frac{1}{x_0}:\dots:\frac{1}{x_n})$ la transformation standard de $\mathbb{P}^n_k$ (de degré $n$), alors $\sigma$ est une involution qui contracte l'hyperplan $H_0$ d'équation $x_0=0$ sur le point $(1:0:\dots:0)$, et qui envoie la valuation associée à $H_0$ sur celle associée au diviseur exceptionnel obtenu en éclatant $(1:0:\dots:0)$. 
En choisissant $p=(1:0:\dots:0)$ (quitte à conjuguer $g$ par un automorphisme de $\mathbb{P}^n_k$), $f=\sigma^{-1} g\sigma\in N$ induit une application birégulière non-triviale de l'hyperplan $H_0$ dans lui-même, correspondant à $g_p\in \Aut(\mathbb{P}(T_p))$. D'après le lemme~\ref{Lem:DefoCremona}, il existe dans $N$ un élément $f_0\in \J_n$, qui préserve la fibration rationnelle $\phi:\mathbb{P}^n_k\dasharrow H_0$ donn\'ee par $(x_0:\dots:x_n)\dasharrow (0:x_1:\dots:x_n)$ et agit sur la base du pinceau comme $f_{|_{H_0}}$, donc de manière non-triviale. 

Montrons maintenant qu'il existe $\beta\in \J^0_n$ tel que $r=\beta^{-1} f_0 \beta (f_0)^{-1}$ soit un élément non-trivial de $N\cap \J_n^0$. Rappelons que $J_n$ est isomorphe au produit semi-direct $J_n^0\rtimes\Bir(k^{n-1})$, et écrivons $f_0=(a,b)$ dans ce produit, avec $a\in J_n^0$ et $b\in \Bir(k^{n-1})$ non trivial par construction. Alors, $r$ s'écrit
\[(\beta^{-1},1)\circ (a,b)\circ (\beta,1)\circ (b^{-1}(a^{-1}),b^{-1})=(\beta^{-1}\cdot a \cdot b(\beta)\cdot a^{-1},1).\]
 Par conséquent, $r$ est un élément de $N\cap J_n^0$, qui est de plus non trivial si et seulement si $\beta\not=a\cdot b(\beta)\cdot a^{-1}$. Si $a$ est l'identité, il suffit de choisir $\beta\in J_n^0$ non fixé par $b$ (par exemple un élément diagonal donné par une fonction de $k(x_1,\dots,x_{n-1})$ qui n'est pas invariante par $b$). Si $a$ n'est pas l'identité, on peut choisir pour $\beta$ un élément de $\PGL(2,k)\subset \PGL(2,k(x_1,\dots,x_{n-1}))$ ne commutant pas avec $a$.

On trouve donc que $N\cap \J^0_n$ est un sous-groupe normal non trivial de $\J_n^0\simeq\PGL(2,k(x_1,\dots,x_{n-1}))$, ce qui implique que $N$ contient $\J_n^1\simeq\PSL(2,k(x_1,\dots,x_{n-1}))$ (voir par exemple \cite[Chapitre II, $\S 2$]{bib:Dieudonne}). De plus, comme $\J^1_n\cap \Aut(\mathbb{P}^n_k)$ est non-trivial et $\Aut(\mathbb{P}^n_k)\simeq \PGL(n+1,k)$ est simple, $N$ contient  $\Aut(\mathbb{P}^n_k)$.
\end{proof}
Rappelons le théorème de Noether-Castelnuovo: le groupe $\Bir(\mathbb{P}^2_k)$ est engendré par $\Aut(\mathbb{P}^2_k)$ et la transformation quadratique standard $(x:y:z)\dasharrow (yz:xz:xy)$ (voir \cite[Chapter V, $\S 5$, Theorem~2, page 100]{bib:Shafa} pour une preuve valable en toute caractéristique). Par conséquent, on a le résultat suivant:
\begin{theo}\label{Thm:Cr2Simple}
Muni de sa topologie, le groupe $\Cr_2(k)=\Bir(\mathbb{P}^2_k)$ est simple.
\end{theo}
\begin{proof}
Suit de la proposition précédente et du fait que $\Bir(\mathbb{P}^2_k)$ soit engendré par $\Aut(\mathbb{P}^2_k)$ et $\J^1_2\simeq \PSL(2,k(x_{1}))$. Démontrons cette dernière partie.  On note $\alpha_{1},\alpha_{2},\beta_{1},\beta_{2}$ les éléments de $\Cr_{2}(k)=\Bir(\mathbb{P}^2_{k})$ suivants (vus ici sur la carte affine $(x_1,x_2)\mapsto (1:x_1:x_2)$ ):
\[\begin{array}{rcclcrccl}
\alpha_{1}\colon& (x_{1},x_{2})&\dasharrow& (x_{1},-\frac{1}{x_{2}}) & \hphantom{espace}&
\beta_{1}\colon& (x_{1},x_{2})&\mapsto& (x_{2},x_{1});\\
\alpha_{2}\colon& (x_{1},x_{2})&\dasharrow& (-\frac{1}{x_{1}},x_{2}) & \hphantom{espace}&
\beta_{2}\colon& (x_{1},x_{2})&\mapsto& (-x_{1},-x_{2}).
\end{array}\]
Alors, $\alpha_{1}$ est un élément de $\J^1_2$ et $\beta_{1},\beta_{2}$ sont deux éléments de $\Aut(\mathbb{P}^2_{k})$. De plus, $\alpha_{2}=\beta_{1}\alpha_{1}\beta_{1}$ et $\alpha_{1}\alpha_{2}\beta_{2}$ est la transformation quadratique standard. Le résultat  se déduit alors du théorème de Noether-Castelnuovo.
\end{proof}
\begin{rema}
La simplicité de $\Bir(\mathbb{P}^2_{k})$, en tant que groupe abstrait, est toujours ouverte. Pour de plus amples résultats dans cette direction, voir \cite{bib:Dan} et \cite{bib:Giz}.
\end{rema}
\section{Connexité de $\Bir(\mathbb{P}^n_k)$}\label{Sec:Connex}
 Comme $\Bir(\mathbb{P}^2_k)$ est engendré par $\Aut(\mathbb{P}^2_k)$ et $\J_2^0$, il est bien connu que le groupe $\Bir(\mathbb{P}^2_k)$ est connexe. En dimension supérieure, il n'existe pas d'analogue au théorème de Noether-Castelnuovo (voir \cite{bib:Pan}) et il ne paraît pas évident à priori de trouver un ensemble adéquat de générateurs. Toutefois, nous pouvons prouver le résultat suivant:

 \begin{theo}\label{Thm:Connexite}
Pour tout $n\geq 1$, le groupe  $\Cr_n(k)=\Bir(\mathbb{P}^n_k)$ est \emph{linéairement connexe} au sens suivant:
pour tous $f,g\in \Cr_n(k)$, il existe un ouvert $U$ de la droite affine sur $k$ contenant $0$ et $1$, et un morphisme $\theta\colon U\to \Cr_n(k)$ tel que $\theta(0)=f$ et $\theta(1)=g$.

En particulier, le groupe $\Cr_n(k)$ est connexe. \end{theo}
 \begin{proof} 
  Si $U\subset \mathbb{A}^1_k$ est un ouvert contenant $0$ et $1$ et que le morphisme $\theta\colon U\to \Cr_n(k)$ satisfait $\theta(0)=f$ et $\theta(1)=g$, on dira que $\theta$ \emph{joint $f$ à $g$}; en notant $U'$ l'ouvert qui est l'image de $U$ par $t\mapsto 1-t$, le morphisme $U'\to \Cr_n(k)$ défini par $t\mapsto \theta(1-t)$ joint $g$ à $f$. De plus, si  $\nu\colon V\to \Cr_n(k)$ joint $g$ à $h$,  le morphisme $U\cap V\to \Cr_n(k)$ défini par $t\mapsto \theta(t)\circ g^{-1}\circ\nu(t)$ joint $f$ à $h$. On en déduit que la relation "$f$ et $g$ sont joignables" est une relation d'équivalence.
 
 Notons $\mathcal{U}_0\subset \Cr_n(k)$ l'ensemble des éléments joignables à l'identité. On observe que $\mathcal{U}_0$ est un sous-groupe normal de $\Cr_{n}(k)$. En effet, si $\theta$ joint $1$ à $f$, alors $t\mapsto \theta(1-t)\circ f^{-1}$ joint $1$ à $f^{-1}$ et si $\nu$ joint $1$ à $g$, alors $t\mapsto \theta(t) \circ \nu(t)$ joint $1$ à $f\circ g$; de plus si $h\in \Cr_n(k)$, $t\mapsto h\circ \theta(t) \circ h^{-1}$ joint $1$ à $h\circ f\circ h^{-1}$.
 
 Montrons maintenant que $\Aut(\mathbb{P}^n_k)\simeq \PGL(n+1,k)$ est contenu dans $\mathcal{U}_0$ (c'est-à-dire qu'il est linéairement connexe). Les éléments de la forme  \[(x_0:x_1:\dots:x_n)\mapsto (x_0+\sum_{i=1}^n a_{0,i}x_i:x_1+\sum_{i=2}^n a_{1,i}x_i:\dots:x_{n-1}+a_{n-1,n}x_n:x_n)\]   sont joignables à l'identité (remplacer tous les $a_{i,j}$ par $t\cdot a_{i,j}$ donne le morphisme souhaité). De même, un élément diagonal 
 \[(x_0:\dots:x_n)\mapsto (a_0x_0:\dots:a_nx_n)\]
 est joignable à l'identité (remplacer $a_i$ par $(a_i-1)t+1$ donne le morphisme souhaité). Comme ces éléments et leurs conjugués engendrent $\Aut(\p^n_k)$, on en déduit que $\Aut(\p^n_k)\subset \mathcal{U}_0$.

 Le même argument  montre que $\J^0_n\simeq \PGL(2,k(x_1,\dots,x_{n-1}))$ est contenu dans $\mathcal{U}_0$.

 Pour $n=1$, $\Bir(\mathbb{P}^1_{k})=\Aut(\mathbb{P}^1_{k})$, qui est linéairement connexe. On  va alors supposer que $n\geq 2$ et que que $\Cr_{n-1}(k)$ est linéairement connexe (en procédant par induction sur $n$). Alors le groupe $J_n$, engendré par $J^0_n$ et $\Cr_{n-1}(k)$, est contenu dans $\mathcal{U}_0$.

 Montrons maintenant que tout élément $g\in\Bir(\mathbb{P}^n_k)$ appartient à $\mathcal{U}_0$, ce qui donnera le résultat souhaité. Quitte à multiplier $g$ par un élément de $\Aut(\mathbb{P}^n_k)\subset \mathcal{U}_{0}$, on peut supposer que $g$ a un point fixe $p$, et que $g$ et son inverse soient régulières en $p$; on supposera de plus après conjugaison que $p=(1:0:\dots:0)$.
 
 Soit $\sigma:(x_0:\dots:x_n)\dasharrow (\frac{1}{x_0}:\dots:\frac{1}{x_n})$ la transformation standard de $\mathbb{P}^n_k$ (de degré $n$), alors $\sigma$ est une involution qui contracte l'hyperplan $H_0$ d'équation $x_0=0$ sur le point $p$. L'élément $f=\sigma^{-1}g\sigma$ induit une application birégulière de $H_0$ dans lui-même. Le lemme~\ref{Lem:DefoCremona} nous donne un élément $f_0\in J_n$ (la dérivée normale de $f$) tel que $f$ et $f_0$ soient joignables; par conséquent $f\in \mathcal{U}_0$. Le groupe $\mathcal{U}_0$ étant normal dans $\Cr_n(k)$, $g$ appartient également à $\mathcal{U}_0$.
\end{proof}

\end{document}